\documentclass[a4paper,12pt,leqno]{amsart}

\usepackage{amsmath}
\usepackage{tikz}
\usepackage{mathdots}
\usepackage{yhmath}
\usepackage{cancel}
\usepackage{color}
\usepackage{siunitx}
\usepackage{array}
\usepackage{multirow}
\usepackage{amssymb}
\usepackage{gensymb}
\usepackage{tabularx}
\usepackage{booktabs}
\usetikzlibrary{fadings}
\usetikzlibrary{patterns}
\usetikzlibrary{shadows.blur}
\usetikzlibrary{shapes}

\RequirePackage{latexsym} \RequirePackage{amsthm}
\usepackage{enumerate}
\usepackage{dsfont}
\usepackage{mathrsfs}
\numberwithin{equation}{section}

\newtheorem{conjecture} {\sc  Conjecture\rm} [section]

\newtheorem{preremark}[conjecture]{Remark}
\newenvironment{remark}%
  {\begin{preremark}}{\end{preremark}}

\newtheorem{predefinition}[conjecture]{Definition}
  {\begin{predefinition}}{\end{predefinition}}

\newtheorem{prelemma}[conjecture]{Lemma}
 {\begin{prelemma}}{\end{prelemma}}
\newtheorem{preproposition}[conjecture]{Proposition}
\newenvironment{proposition}%
 {\begin{preproposition}}{\end{preproposition}}

\newtheorem{precorollary}[conjecture]{Corollary}
\newenvironment{corollary}%
  {\begin{precorollary}}{\end{precorollary}}

\newtheorem{pretheorem}[conjecture]{Theorem}
\newenvironment{theorem}%
 {\begin{pretheorem}}{\end{pretheorem}}

\begin{document}

\title{\bf A Relation of the Allen-Cahn equations and the Euler equations and applications of the equipartition}

\author{Dimitrios Gazoulis}
\address{Department of Mathematics and Applied Mathematics, University of Crete, 70013 Heraklion, Greece, and Institute of Applied and Computational Mathematics, FORTH, 700$\,$13 Heraklion, Crete, Greece}
\email{dgazoulis@math.uoa.gr}
\date{}

\maketitle

\begin{abstract}

We will prove that solutions of the Allen-Cahn equations that satisfy the equipartition of the energy can be transformed into solutions of the Euler equations with constant pressure. As a consequence, we obtain De Giorgi type results, that is, the level sets of entire solutions are hyperplanes. In addition, we obtain some examples of smooth entire solutions of the Euler equations in particular cases. For specific type of initial conditions, some of these solutions can be extended to the Navier-Stokes equations. Also, we will determine the structure of solutions of the Allen-Cahn system in two dimensions that satisfy the equipartition. Finally, we apply the Leray projection on the Allen-Cahn system and provide some explicit entire solutions.

\end{abstract}

$ \\ $

\section{Introduction}

As it is well known, De Giorgi in 1978 \cite{DeGi} suggested a stricking analogy of the Allen Cahn equation ($ \Delta u =f(u) $) with minimal surface theory that led to significant developments in Partial Differential equations and the Calculus of Variations, by stating the following conjecture about bounded solutions on $ {\mathbb{R}}^n : \\ $ 
\textbf{Conjecture:}(De Giorgi) Let $ u \in C^2({\mathbb{R}}^n) $ be a solution to
\begin{align*}
\Delta u - u^3 +u =0
\end{align*}
such that: 1. $ |u|<1 $, $ \; $ 2. $ \dfrac{\partial u}{\partial x_n} >0  \;\: \forall x \in {\mathbb{R}}^n . \\ $ Is it true that all the level sets of $ u $ are hyperplanes, at least for $ n \leq 8 $?

The relationship with the Bernstein problem for minimal graphs is the reason why $ n \leq 8 $ appears in the conjecture.

The conjecture has been proved and the relation of the Allen Cahn with minimal surfaces can be seen via the theory of $ \Gamma $-convergence (see \cite{Savin} for further details). The family of functionals
\begin{align*}
J_{\varepsilon}(u) = \int_{\Omega} ( \frac{\varepsilon}{2} |\nabla u|^2 + \frac{1}{\varepsilon}W(u)) dx \;\; , \; \varepsilon >0
\end{align*}
$ \Gamma$-converges as $ \varepsilon \rightarrow 0 $ to the perimeter functional and the Euler-Lagrange equations are
\begin{align*}
\varepsilon \Delta u - \frac{1}{\varepsilon}W'(u) =0
\end{align*}
therefore one expects that the level sets of the minimizers will minimize the perimeter.

So, one question could be, whether there exists a transformation that transforms the Allen-Cahn equation $ \Delta u = f(u) \;\: (u:\Omega \subset {\mathbb{R}}^n \rightarrow \mathbb{R} $) to the minimal surface equation of one dimension lower (i.e. $ (n-1) $-dimensional minimal surface equation). The answer is positive, at least in dimension 3, by Corollary \ref{Corollary3D} and then we can apply a Bernstein-type theorem for the minimal surface equation (see \cite{CM}, \cite{EW}) to obtain that the level sets of solutions of the Allen-Cahn equation that satisfy the equipartition of the energy are hyperlanes.

For bounded entire solutions of the Allen-Cahn equation that satisfy the equipartition holds a more general result (see Theorem 5.1 in \cite{CGS}), that is, the level sets of entire solutions of the Allen Cahn equations that satisfy the equipartition are hyperplanes. This was already known by Modica and Mortola in 1980 (see final remark in \cite{MM}). In fact, any solution of the Allen-Cahn equation is smooth and satisfies the bound $ | u | \leq 1 $ (see Proposition 1.9 in \cite{Farina}). The point in Corollary \ref{Corollary3D} is that, we can obtain that the level set of solutions are hyperplanes in any open, convex domain with the appropriate boundary conditions, utilizing the result in \cite{EW}. 

As we can see in section 3, the relations between different classes of equations, allow us to obtain some explicit smooth entire solutions for the 2D and 3D Isobaric Euler equations. Those solutions can be extended when the pressure is linear function in the space variables. Some of these solutions have linear dependent components. Thus, if we impose linear dependency of the components of the initial conditions, we can obtain some explicit entire solutions and can be extended to other type of equations. In Appendix B we give some examples of smooth entire solutions of the Navier-Stokes equations with linear dependent components of the initial conditions.

One of the observations in this paper, is to view the equipartition as the Eikonal equation. As stated in Proposition \ref{Eikonal-Euler}, the Eikonal equation can be transformed to the Euler equations with constant pressure (without the divergence free condition). Thus, solutions of the Allen Cahn equations that satisfy the equipartition can be transformed into the Euler equations with constant pressure, and we obtain the divergence free condition from the Allen Cahn equations.

Furthermore, we state this result to the equation $ a(u) \Delta u + b(u) | \nabla u|^2 = c(u) $, under the hypothesis that $ u = \Phi (v) $ for some $ v $ that is also in this class of equation. This hypothesis is quite reasonable since the equation $ a(u) \Delta u + b(u) | \nabla u|^2 = c(u) $ is invariant under such transformations, in the sence that if $ u $ is a solution then $ v= F(u) $ is also in this class of equations.

In the last section, we propose an analogue of the De Giorgi type result for the vector Allen-Cahn equations and we will prove that entire solutions of the Allen-Cahn system in dimension 2 that satisfy the equipartition have such a specific structure. Finally, we apply the Helmholtz-Leray decomposition 
in the Allen-Cahn system and obtain an equation, independent from the potential $ W $. Then we apply the Leray projection (i.e. only the divergence free term from the decomposition) and we can determine explicit entire solutions. In Appendix A, we give some examples of such solutions and compare them to the structure we have obtained from Theorem \ref{TheoremAllenCahnSystem}. One such example, for a particular potential $ W \geq 0 $ with finite number of global minima has the property that $ \lim_{x \rightarrow \pm \infty} u(x,y)= a^{\pm} $, where $ a^{\pm} \in \lbrace W=0 \rbrace $ and $ \lim_{y \rightarrow \pm \infty} u(x,y) = U^{\pm} (x) $ where $ U^{\pm} $ are heteroclinic connections of the system (i.e. $ U^{\pm ''} = W_u(U^{\pm} ) $). If fact, we can have infinitely many such solutions

$ \\ $

\begin{large}
\textbf{Acknowledgments:}
\end{large}
I would like to thank my advisor Professor Nicholas Alikakos for his guidance and inspiration, and also for motivating the study of implications of the equipartition in the Allen-Cahn system. Also, I would like to thank S. Papathanasiou and Professors A. Farina and P. Smyrnelis for their valuable comments on a previous version of this paper, which led to various improvements.

$ \\ $

\section{The Allen-Cahn equation and the equipartition}


\subsection{The equipartition of the energy and the Euler equations}

First, we begin with a transformation that relates the Eikonal equation and the Euler equation (with constant pressure and without the incompressibility condition). Note that $ x_n $ plays the role of the ``time parameter'' and $ x_n \in \mathbb{R} $ instead of $ x_n >0 $. We could choose any of $ x_i \;,\; i=1,...,n \;,\; n \geq 2 $ as a ``time parameter'', supposing the monotonicity condition with respect to $ x_i $.

$ \\ $

\begin{proposition}\label{Eikonal-Euler} Let $ v: \Omega \subset {\mathbb{R}}^n \rightarrow \mathbb{R} $ be a smooth solution of
\begin{equation}\label{Eikonal}
| \nabla v|^2 = G(v) 
\end{equation}
where $ G: \mathbb{R}\rightarrow \mathbb{R} $ is a smooth function and suppose that $ v_{x_n} >0 $. 

Then the vector field $ F=(F_1,...,F_{n-1}) $ where $ F_i = \dfrac{v_{x_i}}{v_{x_n}} \;,\; i=1,...,n-1 $ satisfies the Euler equations
\begin{equation}\label{Euler}
F_{x_n} + F {\nabla}_y F = 0 \;\;,\; y= (x_1,...,x_{n-1})
\end{equation}
\end{proposition}
$ \\ $

\begin{proof}
Differentiating \eqref{Eikonal} over $ x_i $ gives
\begin{equation}\label{Eikonal-EulerEq2}
2 \sum_{j=1}^n v_{x_j} v_{x_j x_i} = G'(v) v_{x_i} \;\;,\; i=1,...,n
\end{equation}

Now we have
\begin{equation}\label{Eikonal-EulerEq4}
F_{i \: x_j} = \frac{v_{x_i x_j} v_{x_n} - v_{x_i} v_{x_n x_j}}{v_{x_n}^2} \;\;,\; j=1,...,n
\end{equation}
\begin{equation}\label{Eikonal-EulerEq5}
\Rightarrow F_j F_{i \: x_j} = \frac{ v_{x_j} v_{x_i x_j} v_{x_n} - v_{x_j} v_{x_i} v_{x_n x_j}}{v_{x_n}^3} 
\end{equation}
Thus, by \eqref{Eikonal-EulerEq4} and \eqref{Eikonal-EulerEq5} (for $ i=1,...,n-1 $), we have
\begin{equation}\label{Eikonal-EulerEq6}
\begin{gathered}
F_{i \: x_n} + \sum_{j=1}^{n-1} F_j F_{i \: x_j} = \\ \frac{v_{x_i x_n} v_{x_n}^2 - v_{x_n} v_{x_i}v_{x_n x_n} + \sum_{j=1}^{n-1} ( v_{x_j} v_{x_i x_j} v_{x_n} - v_{x_j} v_{x_i} v_{x_n x_j} ) }{v_{x_n}^3} = \\ \frac{v_{x_n} \sum_{j=1}^n v_{x_j} v_{x_j x_i} - v_{x_i} \sum_{j=1}^n v_{x_j} v_{x_j x_n} }{v_{x_n}^3}
\end{gathered}
\end{equation}
finally, by \eqref{Eikonal-EulerEq2}, the last equation becomes
\begin{equation}
F_{i \: x_n} + \sum_{j=1}^{n-1} F_j F_{i \: x_j} = \frac{v_{x_n} \frac{G'(v)}{2} v_{x_i} - v_{x_i} \frac{G'(v)}{2} v_{x_n} }{v_{x_n}^3}
=0 
\end{equation}
\begin{equation}\label{Eikonal-EulerEq7}
\Rightarrow F_{i \: x_n} + \sum_{j=1}^{n-1} F_j F_{i \: x_j} =0 \;\;,\; i=1,...,n-1
\end{equation}

\end{proof}

$ \\ $
\begin{remark}\label{RmkEulerEikonal} Note that since $ v_{x_n} >0 $ it holds that $ v( \Omega ) \cap \lbrace G = 0 \rbrace = \emptyset $. Indeed, if $ v(x_0) \in \lbrace G = 0 \rbrace \Rightarrow | \nabla v(x_0) |^2 = 0 $ which contradicts $ v_{x_n} >0 $. So, by setting $ \tilde{v} = P(v) $, where $ P'(v) = \frac{1}{\sqrt{G(v)}} $ we have $ \nabla \tilde{v} = P'(v) \nabla v \Rightarrow | \nabla \tilde{v}|^2 = (P'(v))^2 | \nabla v |^2 \Rightarrow | \nabla \tilde{v}|^2 = 1 . $ Thus $ \tilde{v} $ satisfies $ | \nabla \tilde{v}|^2 =1 $ and $ F_i = \dfrac{v_{x_i}}{v_{x_n}} = \dfrac{{\tilde{v}}_{x_i}}{{\tilde{v}}_{x_n}} . $ 
So, at first, it seems that this transformation can be inverted: $ F_1^2 + ... + F_{n-1}^2 = \dfrac{{\tilde{v}}_{x_1}^2 + ... + {\tilde{v}}_{x_{n-1}}^2}{{\tilde{v}}_{x_n}^2} = \dfrac{1}{{\tilde{v}}_{x_n}^2} -1  \Rightarrow {\tilde{v}}_{x_n} = \dfrac{1}{\sqrt{F_1^2 + ... + F_{n-1}^2 +1 }} $
\begin{equation}\label{Euler-Eikonal}
\Rightarrow \tilde{v} = \int \frac{1}{\sqrt{F_1^2 + ... + F_{n-1}^2 +1 }} dx_n + a(x_1,...,x_{n-1})
\end{equation}
That is, if $ F_i \;,\; i=1,...,n-1 $ satisfies the Euler equations $ F_{x_n} +F {\nabla}_y F =0 $, then $ v $ defined by \eqref{Euler-Eikonal} will satisfy the Eikonal equation. This statement is true for $ n=2 $ (see \cite{AG}). But to generalize for $ n \geq 3 $ it appears that further assumptions are needed. So, the class of solutions of the Euler equations with constant pressure seem to be ``richer'' in some sense than the class of solutions of the Eikonal equation, that is, for every smooth solution of the Eikonal equation, we can obtain a solution of the Euler equation, but not vice versa.
\end{remark}
$ \\ $

\begin{theorem}\label{TheoremAllenCahn-Euler}
Let $ u,v : \Omega \subset {\mathbb{R}}^n \rightarrow \mathbb{R} $ such that $ u_{x_n} >0 $ satisfy the equations 
\begin{equation}\label{Theorem1}
\begin{gathered} a(u)\Delta u + b(u) | \nabla u|^2 = f(u) \\
k(v) \Delta v + l(v) | \nabla v|^2 = g(v)
\end{gathered}
\end{equation}
and suppose that $ u = \Phi(v) $ for some $ \Phi : \mathbb{R} \rightarrow \mathbb{R} \;( \Phi' \neq 0) $ and $ p(t) \neq 0 \;,\; a(t) \neq 0 $, where $ \\ p(t):= k(t) a( \Phi (t)) \Phi ''(t) + k(t) b( \Phi (t)) ( \Phi '(t))^2 - l(t) a( \Phi (t) ) \Phi '(t) \: $. 

Then the vector field $ F=(F_1,...,F_{n-1}) $ defined as $ F_i = \dfrac{u_{x_i}}{u_{x_n}} \:, \\ i=1,..,n-1 $, will satisfy the Euler equations
\begin{equation}\label{Theorem1Euler}
F_{x_n} + F {\nabla}_y F = 0 \;\;\;,\; y=(x_1,...,x_{n-1})
\end{equation}
Also, $div_y F = 0 $ if and only if $ \Phi $ is a solution of the ODE
\begin{align}\label{ODEfordivF=0}
a( \Phi (t)) \Phi '(t) G'(t) + 2 [b( \Phi (t)) ( \Phi '(t))^2 + a( \Phi(t)) \Phi''(t) ] G(t) = 2 f( \Phi (t))
\end{align}
where $ G(t) := \dfrac{k(t) f( \Phi (t)) - g(t) a( \Phi (t)) \Phi '(t)}{p(t)} \;\;\; (p $ as defined above)
\end{theorem} 
$ \\ \\ $

\begin{proof}

We have $ u= \Phi (v) $ and $ \nabla u = \Phi ' (v) \nabla v $ , therefore
\begin{align*}
\begin{gathered}
\Delta u = \Phi'(v) \Delta v + \Phi''(v) |\nabla v|^2 \\
\Rightarrow f(u) - b(u) | \nabla u |^2 = a( \Phi(v)) ( \Phi'(v) \Delta v + \Phi''(v) |\nabla v|^2) \\
\Rightarrow f( \Phi (v)) - b( \Phi (v)) ( \Phi '(v))^2 | \nabla v |^2 = a( \Phi(v)) ( \Phi'(v) \Delta v + \Phi''(v) |\nabla v|^2) \\
\end{gathered}
\end{align*}
\begin{equation}\label{Theorem1proofEq1}
\Rightarrow \Delta v = \dfrac{f( \Phi (v)) - [b( \Phi (v)) ( \Phi '(v))^2 + a( \Phi(v)) \Phi''(v) ] \: | \nabla v |^2 }{a( \Phi (v)) \Phi '(v)} \;\;\; ( a,\Phi' \neq 0)
\end{equation}
since $ u $ is a solution of $ a(u)\Delta u + b(u) | \nabla u|^2 = f(u) . \\ $ 

Now, since $ v $ is also solution of the second equation in \eqref{Theorem1}, we have
\begin{align*}
\begin{gathered}
k(v)(\dfrac{f( \Phi (v)) - [b( \Phi (v)) ( \Phi '(v))^2 + a( \Phi(v)) \Phi''(v) ] \: | \nabla v |^2 }{a( \Phi (v)) \Phi '(v)}) + l(v) | \nabla v |^2 = g(v) \\
\Leftrightarrow p(v) | \nabla v |^2 = k(v)f( \Phi (v)) -a( \Phi(v))  \Phi'(v) g(v)
\end{gathered}
\end{align*}
where $ p(v) = k(v) a( \Phi (v)) \Phi ''(v) + k(v) b( \Phi (v)) ( \Phi '(v))^2 - l(v) a( \Phi (v) ) \Phi '(v) $. 

By hypothesis $ p \neq 0 $, thus
\begin{equation}\label{Theorem1equipartition}
|\nabla v|^2 = G(v)
\end{equation}
where 
\begin{equation}
G(v) = \frac{k(v) f( \Phi (v)) - g(v) a( \Phi (v)) \Phi '(v)}{p(v)}
\end{equation}
Also note that $ F_i = \dfrac{u_{x_i}}{u_{x_n}} = \dfrac{v_{x_i}}{v_{x_n}} . \\ $

So we apply Proposition \ref{Eikonal-Euler} and we obtain that 
\begin{equation}\label{Theorem1Euler2}
\begin{gathered}
F_{i \: x_n} + \sum_{j=1}^{n-1} F_j F_{i \: x_j} =0 \;\;,\; i=1,...,n-1 \\
\Leftrightarrow F_{x_n} + F {\nabla}_y F =0
\end{gathered}
\end{equation}
Now, for the divergence of $ F $:
\begin{align*}
\begin{gathered}
F_{i \: x_i} = \frac{v_{x_i x_i} v_{x_n} -v_{x_i}v_{x_n x_i}}{v_{x_n}^2} \\ \Rightarrow div_y F = \sum_{i=1}^{n-1} F_{i \: x_i} = \frac{\sum_{i=1}^{n-1} v_{x_i x_i} v_{x_n} - \sum_{i=1}^{n-1} v_{x_i} v_{x_n x_i} }{v_{x_n}^2}
\end{gathered}
\end{align*}
\begin{equation}\label{Theorem1div}
\Rightarrow div_y F = \frac{v_{x_n} \Delta v - \frac{1}{2} (|\nabla v|^2)_{x_n}}{v_{x_n}^2}
\end{equation}
Thus, from \eqref{Theorem1proofEq1} and \eqref{Theorem1equipartition} the equation \eqref{Theorem1div} becomes:
\begin{align*}
div_y F = \frac{v_{x_n} \frac{f( \Phi (v)) - [b( \Phi (v)) ( \Phi '(v))^2 + a( \Phi(v)) \Phi''(v) ] \: G(v) }{a( \Phi (v)) \Phi '(v)} - \frac{G'(v)}{2} v_{x_n}}{v_{x_n}^2} 
\end{align*}
Therefore
\begin{align*}
div_y F = 0 \Leftrightarrow a( \Phi (v)) \Phi '(v) G'(v) + 2 [b( \Phi (v)) ( \Phi '(v))^2 + a( \Phi(v)) \Phi''(v) ] G(v) = 2 f( \Phi (v))
\end{align*}
\end{proof}
$ \\ $
\textbf{Notes:} (1) It also holds that solutions of the Allen Cahn equations that satisfy the equipartition also satisfy $ div(\dfrac{\nabla u}{| \nabla u|}) =0 $ has been proved for more general type of equations (see Proposition 4.11 in \cite{DG}). 
$ \\ $
(2) We could see the fact that $ div_y F=0 $, can alternatively be obtained with calculations utilizing the stress-energy tensor (see \cite{AFS} ,p.88), applied in the scalar case. $ \\ $
(3) In the special case of Theorem \ref{TheoremAllenCahn-Euler} where $ u=v \;,\: a(u)=1=l(u) \:,\: b(u)=0=k(u) $ and \eqref{ODEfordivF=0} becomes $ g'(u)=2f(u) $, the equation $ | \nabla u |^2 =g(u) $ is enough to be satisfied only at a single point (see Theorem 2.2 in \cite{CGS}). In fact, it holds that in all dimensions the level sets of bounded entire solutions of the Allen Cahn equations that satisfy the equipartition of the energy are hyperplanes (see the final remark in \cite{MM} or Theorem 5.1 in \cite{CGS}).

$ \\ $

\section{Entire solutions of the Euler equations}

$ \\ $

In this subsection we will determine some smooth entire solutions of the 2D and 3D Euler equations and the pressure being a linear function with respect to the space variables.

We begin by illustrating an analogy for steady solutions of the incompressible Euler equations in two space dimensions and the De Giorgi conjecture.

Let $ u=(u_1,u_2): \mathbb{R}^2 \times (0, + \infty) \rightarrow \mathbb{R}^2 \;,\; u_i =u_i(x,t) \;,\; x=(x_1,x_2) $ be a smooth solution of the Euler equations. The incompressibility condition $ div \: u =0 $ gives that there exists a (unique up to an additive constant) stream function $ \psi(x,t) $  such that
\begin{align*}
u = (- \psi_{x_2}, \psi_{x_1})
\end{align*}
In addition, by Proposition 2.2 in \cite{MB}, a stream function $ \psi $ on a domain $ \Omega \subset \mathbb{R}^2 $ defines a steady solution (i.e. time independent) of the 2D Euler equation on $ \Omega $ if and only if
\begin{align*}
\Delta \psi = F(\psi) \;,\; \textrm{for some function} \;\: F
\end{align*}
So, if $ \psi $ is a bounded, entire solution such that $ \psi_{x_2} \geq 0 $, then by De Giorgi's conjecture (see Theorem 1.1 in \cite{GG}) it holds that
\begin{align*}
\psi(x_1,x_2) = g(ax_1 +bx_2)
\end{align*}
Therefore we raise the following question.
$ \\ $

\textbf{Question:} Let $ u: \mathbb{R}^2 \times (0, + \infty) \rightarrow \mathbb{R}^2 \;,\; (u = u(x,y,t)=(u_1,u_2) )$ be a smooth, bounded entire solution of the Isobaric 2D Euler equations
\begin{align}\label{Euler2DIsobaric}
\begin{cases}
u_{1 \: t} +u_1 u_{1 \: x} + u_2 u_{1 \: y} = 0 \\
u_{2 \: t} +u_1 u_{2 \: x} + u_2 u_{2 \: y} = 0 \\
u_{1 \: x} + u_{2 \: y} =0
\end{cases}
\end{align}

Is it true that then 
\begin{equation}\label{Euler2DSolIsobaric}
u_1 = c_1 \: g(\beta x + \gamma y -( \beta \tilde{c}_1 + \gamma \tilde{c}_2 )t) + \tilde{c}_1 \;\; , \;\: u_2 = c_2 g(\beta x + \gamma y -( \beta \tilde{c}_1 + \gamma \tilde{c}_2 )t) + \tilde{c}_2 \;\;\; ? 
\end{equation}
where $ c_1 \beta + c_2 \gamma =0 \;\;,\; c_1,c_2,\tilde{c}_1,\tilde{c}_2, \beta, \gamma \in \mathbb{R} . \\ \\ $

From the form of solution \eqref{Euler2DSolIsobaric} we can obtain a solution of the 2D Euler equation with pressure being a linear function in respect to the space variables.

Let $ u: \mathbb{R}^3 \rightarrow \mathbb{R}^2 \;,\; (u = u(x,y,t)=(u_1,u_2) )$ is such that
\begin{equation}\label{Euler2DSol}
\begin{gathered}
u_1 = c_1 \: g(\beta x + \gamma y -( \beta \tilde{c}_1 + \gamma \tilde{c}_2 )t) + \lambda A(t) + \tilde{c}_1 \;\; , \;\: u_2 = c_2 g(\beta x + \gamma y -( \beta \tilde{c}_1 + \gamma \tilde{c}_2 )t) + \xi A(t) + \tilde{c}_2 \\ \textrm{and} \;\; p(x,y,t) = -a(t)( \lambda x+ \xi y) +b(t)
\end{gathered}
\end{equation}
where $ A'(t)=a(t) \;,\; a,b: \mathbb{R} \rightarrow \mathbb{R} $ and $ c_1, \tilde{c}_1 , c_2, \tilde{c}_2 ,\beta, \gamma, \lambda , \xi \in \mathbb{R} $ are such that $ c_1 \beta + c_2 \gamma =0 $ and $ \lambda \beta + \xi \gamma =0 . \\ $

Then $ u = (u_1,u_2) $ satisfies
\begin{align}\label{Euler2D}
\begin{cases}
u_{1 \: t} +u_1 u_{1 \: x} + u_2 u_{1 \: y} =-p_x \\
u_{2 \: t} +u_1 u_{2 \: x} + u_2 u_{2 \: y} =-p_y \\
u_{1 \: x} + u_{2 \: y} =0
\end{cases}
\end{align}
$ \\ $

Now we give some examples of smooth entire solutions for the three dimensional Euler equations. If $ u = (u_1,u_2,u_3) : \mathbb{R}^4 \rightarrow \mathbb{R}^3 $ where $ u_i = u_i(x,y,z,t) $ is such that
\begin{equation}\label{Euler3DSol1}
\begin{gathered}
u_1(x,y,z,t) = G(c_1t -y+c_2z) \;\;\;,\;\; u_2(x,y,z,t) = H(c_1t-y+c_2z) -A(t) \\  u_3(x,y,z,t) = \frac{1}{c_2} H(c_1t-y+c_2z) - \frac{1}{c_2} A(t) + C \;\;\textrm{and} \;\; p(x,y,z,t) = a(t) (y + \frac{z}{c_2}) +b(t) \\
\textrm{where} \;\; A'(t)=a(t) \;,\; a,A,G,H : \mathbb{R} \rightarrow \mathbb{R} \;\;,\;\; c_1,c_2 \in \mathbb{R} \;\; \textrm{and} \;\; C = \frac{-c_1}{c_2}
\end{gathered}
\end{equation}
then $ u = (u_1,u_2,u_3) $ is an entire solution of the Euler equations, that is $ u $ satisfies
\begin{equation}\label{Euler3D}
\begin{cases}
u_{1 \: t} +u_1 u_{1 \: x} + u_2 u_{1 \: y}+ u_3 u_{1 \: z} = - p_x \\
u_{2 \: t} +u_1 u_{2 \: x} + u_2 u_{2 \: y} + u_3 u_{2 \: z} = - p_y \\
u_{3 \: t} +u_1 u_{3 \: x} + u_2 u_{3 \: y} + u_3 u_{3 \: z} = -p_z \\
u_{1 \: x} + u_{2 \: y} + u_{3 \: z} =0
\end{cases}
\end{equation}

Note that from symmetry properties of the Euler equations and from \eqref{Euler3DSol1} we can also have the following solution of \eqref{Euler3D}:
\begin{equation}\label{Euler3DSol2}
\begin{gathered}
u_1(x,y,z,t) = \frac{1}{c_2} H(c_1t-z +c_2x) - \frac{1}{c_2} A(t) + C \;\;\;,\;\; u_2(x,y,z,t) = G(c_1t -z +c_2x) \\ u_3(x,y,z,t) = H(c_1t-z +c_2x) -A(t) \;\;\textrm{and} \;\; p(x,y,z,t) = a(t) (z + \frac{x}{c_2}) +b(t) \\
\textrm{where} \;\; A'(t)=a(t) \;,\; a,A,G,H : \mathbb{R} \rightarrow \mathbb{R} \;\;,\;\; c_1,c_2 \in \mathbb{R} \;\; \textrm{and} \;\; C = \frac{-c_1}{c_2}
\end{gathered}
\end{equation}
and also,
\begin{equation}\label{Euler3DSol3}
\begin{gathered}
u_1(x,y,z,t) = H(c_1t-x+c_2y) -A(t) \;,\; u_2(x,y,z,t) = \frac{1}{c_2} H(c_1t-x+c_2y) - \frac{1}{c_2} A(t) + C \\ u_3(x,y,z,t) = G(c_1t -x+c_2y) \;\;\textrm{and} \;\; p(x,y,z,t) = a(t) (x + \frac{y}{c_2}) +b(t) \\
\textrm{where} \;\; A'(t)=a(t) \;,\; a,A,G,H : \mathbb{R} \rightarrow \mathbb{R} \;\;,\;\; c_1,c_2 \in \mathbb{R} \;\; \textrm{and} \;\; C = \frac{-c_1}{c_2}
\end{gathered}
\end{equation}

Finally, another example of smooth entire solution of \eqref{Euler3D} is the following
\begin{equation}\label{Euler3DSol5}
\begin{gathered}
u_1(x,y,z,t) = G([k \tilde{c}_1 + l \tilde{c}_2]t + [k c_1 + l c_2]x -ky -lz)  -A(t) \\
u_2(x,y,z,t) = c_1 u_1(x,y,z,t) + \tilde{c}_1 \;\;\;,\;\;\; u_3(x,y,z,t) = c_2u_1(x,y,z,t) + \tilde{c}_2 \\ \textrm{and} \;\; p(x,y,z,t) = a(t) (x +c_1 y +c_2z) \\
\textrm{where} \;\; A'(t) = a(t) \;\;,\; a,A,G : \mathbb{R} \rightarrow \mathbb{R} \;\; \textrm{and} \;\; c_1,c_2,\tilde{c}_1, \tilde{c}_2,k,l \in \mathbb{R} .
\end{gathered}
\end{equation}
(we can choose $ A $ such that $ A(0)=0 ) \\ $

Therefore we conclude to the following result $ \\ $

\begin{theorem}\label{TheoremEulerEquations}
Let $ u=(u_1,u_2,u_3) \;\:,\; u_i \:,p : \mathbb{R}^3 \times (0,+ \infty ) \rightarrow \mathbb{R} $ and consider the initial value problem
\begin{equation}\label{Euler3DIV}
\begin{cases}
u_t + u \nabla u = - \nabla p \\
\textrm{div} \: u =0 \\
u(x,y,z,0)=g(x,y,z)
\end{cases}
\end{equation}
where $ g = (g_1,g_2,g_3) $ is either of the form
\begin{align}\label{Euler3DIVcond1}
\begin{gathered}
g=(g_1, c_1g_1 + \tilde{c}_1,c_2g_1+ \tilde{c}_2) \;\; \textrm{and} \;\; g_1(x,y,z)= g_1([kc_1 +lc_2]x-ky-lz) \\ c_1,c_2, \tilde{c}_1, \tilde{c}_2 , k,l \in \mathbb{R} \;\;,\; g_1 \;\: \textrm{smooth}
\end{gathered}
\end{align}
or
\begin{align}\label{Euler3DIVcond2}
\begin{gathered}
g=(g_1, g_2,\frac{1}{c_2} g_2 -\frac{c_1}{c_2}) \;\; \textrm{and} \;\; g_1(x,y,z)= G(c_2z-y) \;\;,\;\; g_2(x,y,z) = H(c_2z-y) \\  c_1,c_2, \tilde{c}_1 \in \mathbb{R} \;\:,\;G,H \;\: \textrm{smooth}
\end{gathered}
\end{align}

Then there exists a smooth, globally defined in $ t>0 $, solution of \eqref{Euler3DIV}.

In particular, either $ u $ and $ p $ are given by \eqref{Euler3DSol5} if the initial value $ g $ is of the form \eqref{Euler3DIVcond1} or $ u $ and $ p $ are given by \eqref{Euler3DSol1} if $ g $ is of the form \eqref{Euler3DIVcond2}.
\end{theorem}

$ \\ $

The condition \eqref{Euler3DIVcond2} could be easily modified in order to obtain the solutions given by \eqref{Euler3DSol2} and \eqref{Euler3DSol3}. 

\begin{remark}\label{RmkEulerGeneralDim} Such solutions can be extended to general dimensions, i.e. solutions of \eqref{Euler} and $ n \geq 4 $, together with the divergence free condition and a pressure being a linear function with respect to space variables.
\end{remark}
$ \\ $

\section{Applications}


As we will see now, in dimension 3 we can deduce the Allen Cahn equation (together with the equipartition of the energy) into the minimal surface equation of dimension 2 and then apply Bernstein's result to conclude that the level sets of the solution are hyperplanes. However, the one dimensionality of entire solutions that satisfy the equipartition is a special case of Theorem 5.1 in \cite{CGS} and it was also known by Modica and Mortola (see the final remark in \cite{MM}). However, the result in Corollary \ref{Corollary3D} below holds for any open subset of $ \mathbb{R}^n $, so by imposing the appropriate boundary conditions, utilizing the result in \cite{EW}, we can obtain the result for any convex domain.

$ \\ $

\begin{corollary}\label{Corollary3D} Let $ u \in C^2 (\Omega ; \mathbb{R}) $ be a solution of $ \Delta u = W'(u) $ such that $ u_z >0 $, where $ \Omega \subset \mathbb{R}^3 $ is an open, convex set. If $ u $ satisfies
\begin{align}\label{equipartition2}
\frac{1}{2}|\nabla u|^2 = W(u)
\end{align}
then there exists a function $ \psi $ such that $ \psi_y =- \dfrac{u_x}{u_z} \;,\; \psi_x = \dfrac{u_y}{u_z} $ that satisfies the minimal surface equation
\begin{align*}
\psi_{yy} (\psi_x^2 +1 ) - 2 \psi_x \psi_y \psi_{xy} + \psi_{xx} (\psi_y^2 +1) =0
\end{align*}
In particular, if $ \Omega = {\mathbb{R}}^3 $ or if $ \Omega \subsetneq \mathbb{R}^3 $ and $ u_x=au_z \;,\; u_y =bu_z $ in $ \mathbb{R}^3 \setminus \Omega $, then the level sets of $ u $ are hyperplanes.
\end{corollary}
$ \\ $

\begin{proof}
From Theorem \ref{TheoremAllenCahn-Euler} we have that $ div_{(x,y)} F =0 $, thus there exists some $ \psi = \psi(x,y,z) \: : F_1 = -\psi_y $ and $ F_2 = \psi_x \: .$

As we noted in Remark 2.1, $ u(\Omega) \cap \lbrace W = 0 \rbrace = \emptyset \: $ (by \eqref{equipartition2} and since $ u_z >0 $).

So we set $ v = G(u) $, with $ G'(u) = \dfrac{1}{\sqrt{2W(u)}} $, thus
\begin{align*}
\begin{gathered}
| \nabla v |^2 =1 \;\;\;\; \textrm{and} \;\;\; F_1 = \frac{u_x}{u_z} = \frac{v_x}{v_z} \;,\; F_2 = \frac{u_y}{u_z} = \frac{v_y}{v_z} \\
\Rightarrow F_1^2 +F_2^2 = \frac{1}{v_z^2} -1 \Rightarrow v_z = \frac{1}{\sqrt{F_1^2 +F_2^2 +1}} \\ v_{zx} = \frac{-F_1 F_{1x} - F_2 F_{2x}}{(F_1^2 + F_2^2 +1)^\frac{3}{2}}
\end{gathered}
\end{align*}
and $ v_x = F_1 v_z = \dfrac{F_1}{\sqrt{F_1^2 + F_2^2 +1}} $
\begin{align*}
\Rightarrow v_{xz} = \frac{F_{1z}(F_1^2 +F_2^2 +1) - F_1(F_1 F_{1z} +F_2 F_{2z})}{(F_1^2 + F_2^2 +1)^{\frac{3}{2}}}
\end{align*}

Also, by Proposition \ref{Eikonal-Euler}, $ F $ satisfy
\begin{align*}
\begin{cases} F_{1z} +F_1 F_{1x} + F_2 F_{1y} =0 \\
F_{2z} + F_1 F_{2x} + F_2 F_{2y} =0
\end{cases}
\end{align*}
and therefore, from the fact that $ v_{zx}=v_{xz} $ since $ v \in C^2( \Omega) $, we obtain
\begin{align*}
\begin{gathered}
F_{1z}(F_2^2 +1) - F_1 F_2 F_{2z} + F_1 F_{1x} + F_2 F_{2x} =0 \\
\Rightarrow -F_1 F_2 F_{1x} - F_2^2 F_{1y} - F_{1y} + F_1^2 F_{2x} + F_1 F_2 F_{2y} + F_{2x} =0 \\
\Leftrightarrow \psi_{yy} (\psi_x^2 +1 ) - 2 \psi_x \psi_y \psi_{xy} + \psi_{xx} (\psi_y^2 +1) =0
\end{gathered}
\end{align*}

Finally, if $ \Omega= \mathbb{R}^n $, by Berstein's theorem (see Theorem 1.21 \cite{CM}) $ \psi $ must be a plane (in respect to the variables $ (x,y) $, since $ \psi_{xx}= \psi_{xy} = \psi_{yy}=0 $):
$ {\psi}_x = b(z) $ and $ {\psi}_y =-a(z) $ (for some functions $ a,b: \mathbb{R} \rightarrow \mathbb{R} $)
$ \Rightarrow \psi(x,y,z) = b(z)x -a(z)y + c(z) $. This gives: $ F_1 = -{\psi}_y = a(z) \;\;, \; F_2 = {\psi}_x =b(z) $
\begin{align*}
\begin{gathered}
\Rightarrow \frac{u_x}{u_z} =a(z) \;\;\;\; \textrm{and} \;\;\; \frac{u_y}{u_z} = b(z) \\
\Rightarrow u(x,y,z) = G(s,y) = H(t ,x) \\
\textrm{where} \;\: s=x + \int \frac{1}{a(z)}dz \;\;,\; t =y + \int \frac{1}{b(z)}dz
\end{gathered}
\end{align*}
Now we have
\begin{align*}
\begin{gathered}
u_x = a(z) u_z \Rightarrow H_x = \frac{a}{b} H_t \;\: (H_t \neq 0 \;\: \textrm{since} \;\: u_z >0) \\
\textrm{and} \;\: \frac{1}{2}| \nabla u|^2 =W(u) \Rightarrow \frac{1}{2} [ H_x^2 + H_t^2(1 + \frac{1}{b^2})] = W(H)
\end{gathered}
\end{align*}
Differentiating the last equation with respect to $ y,z $ respectively (and utilizing $ H_x = \frac{a}{b} H_t$), we obtain
\begin{align*}
\begin{cases}
\frac{a}{b} H_t H_{xt} + H_t H_{tt} (1 + \frac{1}{b^2}) = W' H_t \\ \frac{a}{b^2} H_t H_{xt} + H_t H_{tt} ( \frac{1}{b} + \frac{1}{b^3}) - H_t^2 \frac{b'}{b^3} = W' \frac{H_t}{b}
\end{cases}
\end{align*}
\begin{align*}
\Rightarrow - H_t^2 \frac{b'}{b^2} = 0 \Rightarrow b' =0
\end{align*}
thus, $ b=b_0 =constant $. Arguing similarly for $ G= G(s,y) $ we obtain $ a=a_0 = constant. $ Therefore,
\begin{align*}
u(x,y,z) = h(ax+by+z)
\end{align*}
where $ h $ is a solution of the ODE
\begin{align*}
h''(t) = \frac{W'(h(t))}{a^2 +b^2 +1}
\end{align*}
In the case where $ \Omega \subsetneq \mathbb{R}^3 $, we utilize Theorem 1.1 in \cite{EW} to obtain that $ \psi $ is linear in $ \Omega $ and similarly we conclude.

\end{proof}


$ \\ $

Now we will prove an analogue of Theorem 5.1 in \cite{CGS} for subsolutions of the Allen Cahn equation and also, without excluding apriori some potential singularities of the solutions. The observation in Proposition \ref{PropositionAllenCahnEquipartition} below, is to utilize the main result from \cite{CC}.
$ \\ \\ $

\begin{proposition}\label{PropositionAllenCahnEquipartition} Let $ u : \mathbb{R}^n \rightarrow \mathbb{R} $ be a non constant, smooth subsolution of $ \Delta u \leq W'(u) \;, \\ W: \mathbb{R} \rightarrow [0, + \infty) $, except perhaps on a closed set $ S $ of potential singularities with $ \mathcal{H}^1(S) =0 $ and $ \mathbb{R}^n \setminus S $ is connected, such that
\begin{align}\label{PropEquipartition}
\frac{1}{2}|\nabla u|^2 = W(u) 
\end{align}
where $ \mathcal{H}^1 $ is the Hausdorff 1-measure in $ \mathbb{R}^n . $

Then
\begin{align*}
u(x) = g( a \cdot x +b) \;\;\;,\;\; \textrm{for some} \;\: a \in \mathbb{R}^n \;\:,\: |a|=1 \;,\; b \in \mathbb{R}
\end{align*}
and $ g $ is such that $ g'' = W'(g) . $
\end{proposition}

$ \\ $

\begin{proof}
First we see that $ W $ is strictly positive in $ u( \mathbb{R}^n \setminus S) $. Indeed, if there exists $ x_0 \in \mathbb{R}^n \setminus S $ such that $ W(u(x_0)) = 0 $, then $ u $ is a constant by Corollary 3.1 in \cite{AFS} and since $ \mathbb{R}^n \setminus S $ is connected.

So let $ v = G(u) $, where $ G'(u) = \dfrac{1}{\sqrt{2W(u)}} $, then
\begin{align*}
| \nabla v |^2 = (G'(u))^2 | \nabla u |^2 = 1 \;\;\;,\;\; \textrm{on} \;\: \mathbb{R}^n \setminus S
\end{align*}
so $ v $ is a smooth solution of the Eikonal equation except perhaps of a closed set $ S $ of potential singularities with $ \mathcal{H}^1 (S) =0 $. Thus from the result of \cite{CC}, we have that $ v=  a \cdot x +b \;\;,\; a \in \mathbb{R}^n \;,\: |a|=1 \;, \; b \in \mathbb{R} $ or $ v = | x-x_0 | +c $ for some $ x_0 \in \mathbb{R}^n \;\;,\; c \in \mathbb{R} . \\ $

Therefore,
\begin{align*}
u(x) = G^{-1} (a \cdot x +b) \;\;\: \textrm{or} \;\: u(x) = G^{-1} (| x-x_0 | +c)
\end{align*}
where $ G : \mathbb{R} \rightarrow \mathbb{R} $, such that $ G' = \dfrac{1}{\sqrt{2W}} $.

If $ u = G^{-1}(d +c) $ where $ d(x)= |x-x_0|$, then
\begin{align*}
\Delta u = (G^{-1})''(d +c) + \frac{n-1}{d} (G^{-1})'(d+c) \leq W'(u) = W'(G^{-1}(d+c))
\end{align*}
and also,
\begin{align*}
| \nabla u |^2 = (G^{-1 '}(d+c))^2 = 2 W(u) \Rightarrow (G^{-1})'(d+c) = \sqrt{2W(G^{-1}(d+c))}
\end{align*}
and thus, $ (G^{-1})'' = W'(G^{-1}) $, so we obtain
\begin{align*}
(G^{-1})''(d +c) + \frac{n-1}{d} (G^{-1})'(d+c) \leq (G^{-1})''(d+c) \\ \Rightarrow (G^{-1})'(d+c) = 0 \Rightarrow \sqrt{2W(G^{-1}(d+c))} =0
\end{align*}
which contradicts the fact that $ W $ is strictly positive in $ u(\mathbb{R}^n) $.

Therefore $ u(x) = g (a \cdot x +b) $ where $ g = G^{-1} $.

\end{proof}

\begin{remark}\label{RmkSymmetrySolAC} 
\textbf{(1)} In Proposition 2 above, radially symmetric solutions are excluded as we see in the proof, but as it is well known (see \cite{GNN}) if $ f $ is smooth and $ u \in C^2 ( \overline{ \Omega}) $ is a positive solution of $ - \Delta u = f(u) $ for $ x \in B_1 \subset \mathbb{R}^n $ that vanishes on $ \partial B_1 $, it holds that then $ u $ is radially symmetric. So radially symmetric solutions of the Allen-Cahn equations are incompatible with the equipartition even if we do not exclude apriori singularities. 
$ \\ $
\textbf{(2)} Note that, in Theorem \ref{TheoremAllenCahn-Euler}, if $ u,v $ are smooth entire solutions, by \eqref{Theorem1equipartition} in the proof and the monotonicity $ u_{x_n}>0 $, arguing as in the proof of Proposition \ref{PropositionAllenCahnEquipartition} above we can conclude that $ u, \: v $ are one dimensional and the radially symmetric solutions are also excluded in this case.
\end{remark}


$ \\ $

\section{The Allen Cahn system}

\subsection{Applications of the Equipartition}

$ \\ $

We begin by proposing a De Giorgi like result for the Allen Cahn systems for solutions that satisfy the equipartition of the energy or as an analogy of \cite{CGS} in the vector case. First, the property that the level sets of a solution are hyperplanes can be expressed equivalently as $ \dfrac{u_{x_i}}{u_{x_n}} = c_i \;,\;\: i =1,...,n-1 \;\: (u: \mathbb{R}^n \rightarrow \mathbb{R} \;\: , \; u_{x_n}>0 $), that is, if we consider $ v_i = \dfrac{u_{x_i}}{u_{x_n}} \;\:,\; i=1,...,n,\;\: v_i : \mathbb{R}^n \rightarrow \mathbb{R} $, then
\begin{align*}
 v_i = c_i \;\:,\; i=1,...,n-1 \Leftrightarrow rank( \nabla v_i) <1 \;\:,\; i=1,...,n-1
\end{align*}
We can see the above statement as follows,
$ \\ $ If $  v_i = c_i \;\:,\; i=1,...,n-1 $ then $ \nabla v_i =0 \Rightarrow rank( \nabla v_i) =0 <1 \\ $ Conversely, if $ rank( \nabla v_i) <1 $, that is $ rank( \nabla v_i) =0 $ since $ v_i : \mathbb{R}^n \rightarrow \mathbb{R} $, we have by Sard's theorem that $ \mathcal{L}^1(v_i(\mathbb{R}^n)) = 0 \;\:,\; i=1,...,n-1 $ (see for example \cite{RN}) where $ \mathcal{L}^1 $ is the lebesgue measure on $ \mathbb{R}$. Thus, $\mathcal{L}^1(v_i(\mathbb{R}^n)) =0 \Rightarrow v_i =c_i $ (constant) $i=1,...,n-1. \\ $

Now, we can generalize the above to the vector case as follows: $ \\ $

Let $ u : \mathbb{R}^n \rightarrow \mathbb{R}^m \;\:, \; u = (u_1,...,u_m) \;\:,\; u_i = u_i (x_1,...,x_n) $, we consider the functions
\begin{align*}
v_{ij} = \frac{u_{i \: x_j}}{u_{i \: x_n}} \;\:, \; i =1,...,m \;\:,\; j =1,...,n-1
\end{align*}
and $ \tilde{v}^k = (v_{1k},...,v_{mk}) \;\:,\; \tilde{v}^k : \mathbb{R}^n \rightarrow \mathbb{R}^m \;\:,\; k =1,...,n-1 $ and $ \nabla \tilde{v}^k : \mathbb{R}^n \rightarrow \mathbb{R}^{m \times n} $. Thus, if $ u $ is a solution of the Allen Cahn system, we could ask (under appropriate assumptions) whether $ rank ( \nabla \tilde{v}^k) < \min \lbrace n,m \rbrace = \mu $ (and by Sard's  Theorem we would have that $ \mathcal{L}^{\mu} ( \tilde{v}^k (\mathbb{R}^n)) =0 $, where $ \mathcal{L}^{\mu} $ is the Lebesgue measure in $ \mathbb{R}^{\mu}). \\ $

Apart from $ u $ being a solution of the Allen Cahn system (and $ u_{i \: x_n}>0) $) we should need further assumptions, as in the scalar case. The geometric analog in the vector case is far more complicated than in the scalar case. In particular, there is a relationship with minimizing partitions. However, one possible assumption would be that $ u $ also satisfies the equipartition, i.e. $ \dfrac{1}{2} | \nabla u |^2 = W(u) $. We will now prove that the above is true, at least for $ n=m=2 $, that is, if $ \tilde{v} = (v_1,v_2) \;,\; v_i = \dfrac{u_{i \: x}}{u_{i \:y}} $ and $ u=(u_1,u_2) $ is a solution of the Allen-Cahn system that satisfy the equipartition, then $ rank(\nabla v) <2 $. In fact, we can obtain a quite stronger result about the structure of solutions in two dimensions, as stated in Theorem \ref{TheoremAllenCahnSystem} that follows.

$ \\ \\ $

\begin{theorem}\label{TheoremAllenCahnSystem}
Let $ u : \mathbb{R}^2 \rightarrow \mathbb{R}^2 $ be a smooth solution of
\begin{equation}\label{Theorem2AllenCahn}
\Delta u = W_u(u)
\end{equation}
with $ u_{iy}>0 \;,\; i =1,2 $ and $ W: \mathbb{R}^2 \rightarrow [0, +\infty) $ smooth. $ \\ $
If $ u $ satisfies
\begin{equation}\label{Theorem2equipartition}
\frac{1}{2} | \nabla u |^2 = W(u)
\end{equation}
Then
\begin{equation}\label{Theorem2Sol1}
\textrm{either} \;\;\;\: u(x,y) = (U_1 (c_1 x +y), U_2(c_2 x+y)) \;\;\;,\;\: \textrm{where} \;\: U_i'' = \frac{W_{U_i}(U_1,U_2)}{c_i^2 +1} \;\: i=1,2
\end{equation}
\begin{equation}\label{Theorem2Sol2}
\textrm{or} \;\;\;\;\;\;\;\; 
\begin{cases}  h(\dfrac{u_{1x}}{u_{1y}} , \dfrac{u_{2x}}{u_{2y}}) =0 \;\;\;, \\ \textrm{and}  \;\;\; u_{2y}^2 h_{v_1} - u_{1y}^2 h_{v_2} =0
\end{cases}
\end{equation}
for some $ \;\: h: \mathbb{R}^2 \rightarrow \mathbb{R} . $

In particular, $ \mathcal{L}^2(v(\mathbb{R}^2 ))=0 $, where $ v= (\dfrac{u_{1x}}{u_{1y}} , \dfrac{u_{2x}}{u_{2y}}). $
\end{theorem}
$ \\ $

\begin{proof}
We differentiate \eqref{Theorem2equipartition} with respect to $ x \;,\: y $
\begin{equation}\label{Theorem2proofeq1}
\begin{cases} u_{1x} u_{1xx} + u_{1y} u_{1yx} + u_{2x} u_{2xx} + u_{2y} u_{2yx} = W_{u_1} u_{1x} + W_{u_2} u_{2x} \\
u_{1x} u_{1xy} + u_{1y} u_{1yy} + u_{2x} u_{2xy} + u_{2y} u_{2yy} = W_{u_1} u_{1y} +W_{u_2} u_{2y}
\end{cases}
\end{equation}
and utilizing \eqref{Theorem2AllenCahn} we get
\begin{equation}\label{Theorem2proofeq2}
\begin{cases} u_{1x} u_{1xx} + u_{1y} u_{1yx} + u_{2x} u_{2xx} + u_{2y} u_{2yx} = u_{1x} \Delta u_1 + u_{2x} \Delta u_2 \\
u_{1x} u_{1xy} + u_{1y} u_{1yy} + u_{2x} u_{2xy} + u_{2y} u_{2yy} = u_{1y} \Delta u_1 + u_{2y} \Delta u_2
\end{cases}
\end{equation}
\begin{equation}\label{Theorem2proofeq3}
\Leftrightarrow \begin{cases}
u_{1y} u_{1yx} + u_{2y} u_{2yx} = u_{1x}u_{1yy} +u_{2x} u_{2yy} \\
u_{1x} u_{1xy} + u_{2x} u_{2xy} = u_{1y} u_{1xx} + u_{2y} u_{2xx}
\end{cases}
\end{equation}
Now we define $ v_i := \dfrac{u_{ix}}{u_{iy}} \;,\; i=1,2 $ and by the second equation in \eqref{Theorem2proofeq3} we have
\begin{equation}\label{Theorem2proofeq4}
\begin{gathered}
v_{1x}= \frac{u_{1xx}u_{1y}-u_{1x}u_{1yx}}{u_{1y}^2} = \frac{u_{2x}u_{2xy}-u_{2y}u_{2xx}}{u_{1y}^2} = -\frac{u_{2y}^2}{u_{1y}^2} v_{2x} \\
\Leftrightarrow u_{1y}^2 v_{1x} + u_{2y}^2 v_{2x} = 0
\end{gathered}
\end{equation}
similarly by the first equation in \eqref{Theorem2proofeq3} we have
\begin{equation}\label{Theorem2proofeq5}
u_{1y}^2 v_{1y} + u_{2y}^2 v_{2y} = 0
\end{equation}
From \eqref{Theorem2proofeq4}, \eqref{Theorem2proofeq5} and the assumption $ u_{iy}>0 \;,\; i=1,2 $ we obtain that
\begin{equation}
v_{1x} v_{2y} - v_{1y}v_{2x} =0 \Leftrightarrow det( \nabla v) =0 \;\;,\; \forall \; (x,y) \in \mathbb{R}^2 
\end{equation}
Since $ det( \nabla v) = 0 $, we have that $ rank ( \nabla v)<2 $ and by Sard's Theorem (see for example \cite{RN}, p. 20) we have that $ \mathcal{L}^2(v(\mathbb{R}^2 ))=0 $. By Theorem 1.4.14 in \cite{RN}, since $ rank ( \nabla v)<2 $, we have that $ v_1,v_2 $ are functionally dependent, that is, there exists a smooth function $ h: \mathbb{R}^2 \rightarrow \mathbb{R} $ such that
\begin{equation}\label{Theorem2proofeq6}
h(v_1 ,v_2) = 0 \Leftrightarrow h(\frac{u_{1x}}{u_{1y}} , \frac{u_{2x}}{u_{2y}}) =0 \;\;,\; \forall \; (x,y) \in \mathbb{R}^2 
\end{equation}
Thus we have
\begin{equation}\label{Theorem2proofeq7}
h_{v_1}v_{1x} + h_{v_2}v_{2x} =0 \;\; \textrm{and} \;\: h_{v_1}v_{1y} + h_{v_2} v_{2y} =0
\end{equation}
so, together with \eqref{Theorem2proofeq4}, \eqref{Theorem2proofeq5} we get
\begin{equation}\label{Theorem2proofeq8}
(u_{1y}^2 h_{v_2} - u_{2y}^2 h_{v_1}) v_{2x} =0 \;\; \textrm{and} \;\: (u_{1y}^2 h_{v_2} - u_{2y}^2 h_{v_1}) v_{2y} =0 
\end{equation}
which gives
\begin{equation}\label{Theorem2proofeq9}
\begin{gathered}
v_{2x} = 0 \;\; \textrm{and} \;\; v_{2y} = 0 \\
\textrm{or} \;\; u_{1y}^2 h_{v_2} - u_{2y}^2 h_{v_1}=0
\end{gathered}
\end{equation}
in the first case we also have
\begin{equation}\label{Theorem2proofeq10}
v_{1x} = 0 \;\; \textrm{and} \;\; v_{1y} = 0 
\end{equation}
and therefore
\begin{equation}\label{Theorem2proofeq11}
\begin{gathered}
\frac{u_{1x}}{u_{1y}} = c_1  \;\; \textrm{and} \;\; \frac{u_{2x}}{u_{2y}} =c_2
\Rightarrow u_1 (x,y) = U_1 (c_1 x +y) \;\; \textrm{and} \;\; u_2(x,y) = U_2 (c_2 x +y)
\end{gathered}
\end{equation}
where
\begin{align*}
U_i'' = \frac{W_{U_i}(U_1,U_2)}{c_i^2 +1} \;\: i=1,2.
\end{align*}
In the second case we see that both equations of \eqref{Theorem2Sol2} are satisfied.
\end{proof}

$ \\ \\ $
\textbf{Note:} If $ W(u_1,u_2) = W_1(u_1) +W_2(u_2) $, then \eqref{Theorem2AllenCahn} becomes
\begin{align*}
\Delta u_i = W_i'(u_i) \;\;,\; i=1,2
\end{align*}
so, by analogy with the scalar case we should suppose $ u_{iy} >0 $ as we see in Theorem \ref{TheoremAllenCahnSystem} above.
$ \\ \\ $

\subsection{The Leray projection on the Allen-Cahn system}

$ \\ $

We begin with a calculation with which we will obtain an equation independent of the potential $ W . \\ $

Let $ u : \mathbb{R}^2 \rightarrow \mathbb{R}^2 $ be a smooth solution of the system
\begin{equation}\label{AllenCahnSystem}
\Delta u = W_u(u) \Leftrightarrow \begin{cases} \Delta u_1 = W_{u_1}(u_1,u_2) \\ \Delta u_2 = W_{u_2} (u_1,u_2) \end{cases}
\end{equation}
where $ W : \mathbb{R}^2 \rightarrow \mathbb{R} $, a $ C^2 $ potential.

From \eqref{AllenCahnSystem}, differentiating over $ x \;,\:y $ we obtain
\begin{equation}\label{DifferAllenCahnSystem1}
\begin{cases}
\Delta u_{1y} = W_{u_1 u_1} u_{1y} + W_{u_1 u_2} u_{2y} \\
\Delta u_{1x} = W_{u_1 u_1} u_{1x} + W_{u_1 u_2} u_{2x} \\
\Delta u_{2y} = W_{u_2 u_1} u_{1y} + W_{u_2 u_2} u_{2y} \\
\Delta u_{2x} = W_{u_2 u_1} u_{1x} + W_{u_2 u_2} u_{2x}
\end{cases}
\end{equation}
and therefore
\begin{equation}\label{DifferAllenCahnSystem2}
u_{1x} \Delta u_{1y} + u_{2x} \Delta u_{2y} = W_{u_1 u_1} u_{1y} u_{1x} + W_{u_1 u_2} ( u_{1x} u_{2y} + u_{1y} u_{2x} ) + W_{u_2 u_2} u_{2y} u_{2x}
\end{equation}
thus we have
\begin{equation}\label{DifferAllenCahnSystem3}
u_{1x} \Delta u_{1y} + u_{2x} \Delta u_{2y} = u_{1y} \Delta u_{1x} + u_{2y} \Delta u_{2x}
\end{equation}

Now we will apply the Helmholtz-Leray decomposition, that resolves a vector field $ u $ in $ \mathbb{R}^n \;\: (n = 2, 3) $ into the sum of a gradient and a curl vector. Regardless of any boundary conditions, for a given vector field $ u $ can be decomposed in the form
\begin{align*}
u = \nabla \phi + \tilde{\sigma} = ( \phi_x + \tilde{\sigma}_1 , \phi_y + \tilde{\sigma}_2)
\end{align*}
where $ div \: \tilde{\sigma} =0 \Leftrightarrow \tilde{\sigma}_{1 x} + \tilde{\sigma}_{2 y} =0  $ since we are in two dimensions, and thus $ \tilde{\sigma}_1 = - \sigma_y \;,\; \tilde{\sigma}_2 = \sigma_x $. So, we have that
\begin{align*}
u = ( \phi_x - \sigma_y , \phi_y +  \sigma_x)
\end{align*}
for some $ \phi , \sigma : \mathbb{R}^2 \rightarrow \mathbb{R} $.

Utilizing now this decomposition of $ u $, we obtain
\begin{equation}\label{Helmholtz-LerayDecomEq}
\begin{gathered}
(\phi_{xx} - \sigma_{yx}) \Delta (\phi_{xy} - \sigma_{yy}) + ( \phi_{yx} + \sigma_{xx}) \Delta (\phi_{yy} + \sigma_{xy}) \\ = (\phi_{xy} - \sigma_{yy}) \Delta(\phi_{xx} - \sigma_{yx}) + ( \phi_{yy} + \sigma_{xy}) \Delta( \phi_{yx} + \sigma_{xx})
\end{gathered}
\end{equation}
Thus, if in particular we apply the Leray projection, $ v = \mathbb{P}(u) $, we have that $ v = \tilde{\sigma} $, that is, $ v = ( - \sigma_y , \sigma_x) $. So, from \eqref{Helmholtz-LerayDecomEq} we have
\begin{equation}\label{LerayPrEq1}
\begin{gathered}
\sigma_{yx} \Delta \sigma_{yy} + \sigma_{xx} \Delta \sigma_{xy} = \sigma_{yy} \Delta \sigma_{yx} +  \sigma_{xy} \Delta \sigma_{xx} \\ \Leftrightarrow (\sigma_{xx} - \sigma_{yy}) \Delta \sigma_{xy} = \sigma_{xy} \Delta ( \sigma_{xx} - \sigma_{yy})
\end{gathered}
\end{equation}

Note that a class of solutions to \eqref{LerayPrEq1} is $ \sigma $ that satisfy
\begin{equation}\label{LerayPrEq2}
c_1 \sigma_{xy} = c_2 ( \sigma_{xx} - \sigma_{yy})
\end{equation}
and we can solve explicitly in $ \mathbb{R}^2 $,
\begin{equation}\label{LerayPrSol}
\begin{gathered}
\sigma(x,y) = A(x) + B(y) \;\;\;,\;\; \textrm{if} \;\: c_2=0 \\
\sigma(x,y) = F(cx +y) + G(x - cy) \;\;\;,\;\; \textrm{where} \;\: c = \dfrac{c_1 + \sqrt{c_1^2 + 4c_2^2}}{2c_2} \;,\;\: \textrm{if} \;\: c_2 \neq 0
\end{gathered}
\end{equation}
for arbitrary functions $ A,B,F,G : \mathbb{R} \rightarrow \mathbb{R} $.

In the first case, the Leray projection of the solution is of the form
\begin{align}\label{LerayPrSol1}
v = \mathbb{P}(u) = ( b(y) , a(x))
\end{align}
and in the second case
\begin{align}\label{LerayPrSol2}
v = \mathbb{P}(u) = ( c g(x -cy) - f(cx+y), g(x -cy) + cf(cx +y))
\end{align}

Similarly, if we take the projection to the space of gradients, we have $ \tilde{v} = ( \phi_x , \phi_y) $ that will also satisfy
\begin{equation}\label{GradProjEq}
(\phi_{xx} - \phi_{yy}) \Delta \phi_{xy} = \phi_{xy} \Delta ( \phi_{xx} - \phi_{yy})
\end{equation}
so again, the projection to the space of gradients of the solution will be of the form
\begin{equation}\label{GradProjSol}
\begin{gathered}
\textrm{either} \;\: \tilde{u}(x,y) = (\tilde{a}(x),\tilde{b}(y)) \\
\textrm{or} \;\: \tilde{u}(x,y) = ( c \tilde{f}(cx+y) +\tilde{g}(x-cy), \tilde{f}(cx+y) - c \tilde{g}(x-cy))
\end{gathered}
\end{equation}

Therefore, if we determine a class of potentials $ W $, such that the solutions (or some solutions) are invariant under the Leray projection (or the projection to the space of gradients), we can obtain explicit solutions of the form \eqref{LerayPrSol} or \eqref{GradProjSol}. In the Appendix we give such examples.

$ \\ \\ \\ \\ $

\begin{appendix}
\begin{LARGE}
\textbf{Appendix}
\end{LARGE}
\section{Some examples of entire solutions of the Allen-Cahn system}
\label{sec:appendix}

$ \\ $

We note that solutions of the form \eqref{LerayPrSol} and \eqref{GradProjSol} are equivalent in the special case that \eqref{LerayPrEq2} is satisfied. So in the class of solutions of \eqref{LerayPrEq2} the Leray projection is, in some sense equivalent with the projection to the space of gradients.
Suppose now that $ u = \nabla \phi $ for some $ \phi : \mathbb{R}^2 \rightarrow \mathbb{R} $, that is, a solution of the Allen-Cahn system remains invariant under the projection to the space of gradients. Then, as \eqref{GradProjEq} we have
\begin{equation}
\phi_{xy} \Delta (\phi_{xx} - \phi_{yy}) = (\phi_{xx} - \phi_{yy}) \Delta \phi_{xy}
\end{equation}

So a simple solution to (A.1) is
\begin{equation}
\phi_{xx} - \phi_{yy} = 0 \Rightarrow \phi (x,y) = F(x+y) + G(x-y)
\end{equation}
and $ u(x,y) = ( \phi_x , \phi_y) $, so in this case $ u $ has the form
\begin{equation}
u(x,y) = (f(x+y) + g(x-y),f(x+y) - g(x-y))
\end{equation}
for some $ f,g : \mathbb{R} \rightarrow \mathbb{R} . \\ $

If $ u $ has the form (A.3), we can see that it also satisfies the equipartition. Indeed, \eqref{AllenCahnSystem} becomes
\begin{equation}
\begin{cases}
2f'' +2g'' = W_{u_1} \\
2f'' - 2g'' = W_{u_2}
\end{cases}
\Rightarrow
\begin{cases}
2(f'' +g'')(f' +g') = W_{u_1}(f' +g') \\
2(f'' - g'')(f' -g') = W_{u_2}(f' -g')
\end{cases}
\end{equation}
\begin{equation}
\begin{gathered}
\Rightarrow 4f''f' +4g''g' = W_{u_1}(f' +g') + W_{u_2}(f' -g') \\
\Rightarrow 2 (f')^2 +2(g')^2 = W(f+g,f-g) +c
\end{gathered}
\end{equation}
$ \\ $
and the equipartition can be written as
\begin{equation}
\begin{gathered}
\frac{1}{2} | \nabla u |^2 = W(u) \\ \Leftrightarrow 2 (f'(x+y))^2 + 2 (g'(x-y))^2 = W(f(x+y)+g(x-y),f(x+y)-g(x-y))
\end{gathered}
\end{equation}
$ \\ $
(the system (A.1) remains equivalent if we add a constant to the potential)
$ \\ $

First we note that solutions of the form (A.3) satisfy \eqref{Theorem2Sol2} in Theorem \ref{TheoremAllenCahnSystem}. Indeed, if $ u $ is of the form (A.3),
\begin{equation}
\begin{gathered}
u_1 = f(x+y) +g(x-y) \;\: \textrm{and} \;\: u_2 = f(x+y) -g(x-y) \\
\Rightarrow \frac{u_{1x}}{u_{1y}} = \frac{f'+g'}{f'-g'} =v_1 \;\;\: \textrm{and} \;\;\: \frac{u_{2x}}{u_{2y}} = \frac{f'-g'}{f'+g'} =v_2
\end{gathered}
\end{equation}
so the function $ h :\mathbb{R}^2 \rightarrow \mathbb{R} $ in \eqref{Theorem2Sol2} is $ h(s,t) = st -1 $. Also,
\begin{equation}
\frac{u_{1y}^2}{u_{2y}^2} = \frac{(f' -g')^2}{(f'+g')^2} \;\;\: \textrm{and} \;\;\: \dfrac{h_{v_1}}{h_{v_2}} = \dfrac{v_2}{v_1} = \dfrac{(f' -g')^2}{(f'+g')^2} = \frac{u_{1y}^2}{u_{2y}^2}
\end{equation}

Now we will see some examples of solutions to the Allen Cahn system that are not in the form \eqref{Theorem2Sol1} in Theorem \ref{TheoremAllenCahnSystem} (which are more similar to the ones in the scalar case). Some of the examples of such solutions are in the form (A.3) and for all solutions in this form the function $ h $ in \eqref{Theorem2Sol2} is, as mentioned above, $ h(s,t) = st -1 . \\ \\ $

\textbf{Example (1)} If $ W(u_1,u_2) = u_1 u_2 $, then 
\begin{align*}
u(x,y) = (\textrm{cosh} (\dfrac{x+y}{\sqrt{2}}) + \textrm{sin}(\dfrac{x-y}{\sqrt{2}}), \textrm{cosh}(\dfrac{x+y}{\sqrt{2}}) - \textrm{sin}(\dfrac{x-y}{\sqrt{2}})) 
\end{align*}
where $ \textrm{cosh} (t) = \dfrac{e^{t} + e^{-t}}{2} $, is a solution of $ \Delta u = W_u(u) $ that satisfies the equipartition and is of the form \eqref{Theorem2Sol2}. A more general solution is
\begin{align*}
u(x,y) = (c_1 e^{a_1x +b_1y} +c_2 e^{a_2x +b_2y} +c_3 \textrm{sin}(a_3x +b_3y) + c_4 \textrm{cos}(a_4x +b_4y), \\ c_1 e^{a_1x +b_1y} +c_2 e^{a_2x +b_2y} -c_3 \textrm{sin}(a_3x +b_3y) - c_4 \textrm{cos}(a_4x +b_4y)) \\ \textrm{where} \;\;\; a_i^2+b_i^2=1 \;,\;\; i=1,2,3,4 \;\:,\; c_i \in \mathbb{R} \;\:  i=1,2,3,4.
\end{align*}
However, not all solutions in this form satisfy the equipartition. In this example the zero set of the potential is $ \lbrace W= 0 \rbrace = \lbrace u_1=0 \rbrace \cup \lbrace u_2=0 \rbrace $. Such potentials $ W $ belong in a class of potentials that have been thoroughly studied in \cite{CL}.
$ \\ \\ $

\textbf{Example (2)} If $ W(u_1,u_2)= \dfrac{[(u_1+u_2)^2-4]^2+[(u_1-u_2)^2-4]^2}{16} $, then
\begin{align}
u(x,y) = (\textrm{tanh}(\dfrac{x+y}{\sqrt{2}}) + \textrm{tanh}(\dfrac{x-y}{\sqrt{2}}), \textrm{tanh}(\dfrac{x+y}{\sqrt{2}}) - \textrm{tanh}(\dfrac{x-y}{\sqrt{2}}))
\end{align}
is a solution of $ \Delta u = W_u(u) $ that satisfies the equipartition (and is of the form \eqref{Theorem2Sol2} and $ h(s,t) = st-1 $). In addition, $ u $ above connects all four phases of the potential $ W $ at infinity, that is
\begin{align*}
\lim_{x \rightarrow \pm \infty} u(x,y) = ( \pm 2 ,0) \;\;\; \textrm{and} \;\;\; \lim_{y \rightarrow \pm \infty} u(x,y) = ( 0, \pm 2 )
\end{align*} 
$ \lbrace W=0 \rbrace = \lbrace (2 ,0), (-2 ,0) , (0, 2 ) , (0, -2) \rbrace . $ 

This solution is a saddle solution (see \cite{Fusco}) and is invariant under rotations of $ \frac{\pi}{2} $ angle (i.e. $ u( \omega (x,y)) = \omega u(x,y) $, where $ \omega $ is the $ \frac{\pi}{2} $-rotation matrix.

Also, another solution of $ \Delta u = W_u(u) $ for such potential is
\begin{align}
u(x,y)=  (\textrm{tanh}x + \textrm{tanh}(\dfrac{x+y}{\sqrt{2}}) , \textrm{tanh}x - \textrm{tanh}(\dfrac{x+y}{\sqrt{2}}))
\end{align} 
for this solution the function $ h $ in \eqref{Theorem2Sol2} is $ h(s,t)=s+t -2 $ but $ u $ in (A.10) does not satisfy the equipartition. Thus, the class of solutions of the Allen-Cahn system that are of the form \eqref{Theorem2Sol2} in Theorem \ref{TheoremAllenCahnSystem}, is more general than that of solutions to the Allen-Cahn system that satisfy the equipartition.
Note that $ u $ in (A.10) has the property that
\begin{equation}
\lim_{x \rightarrow \pm \infty} u(x,y) = ( \pm 2 ,0) \;\;\; \textrm{and} \;\;\; \lim_{y \rightarrow \pm \infty} u(x,y) = (\textrm{tanh}x \pm 1, \textrm{tanh}x \mp 1)
\end{equation}
and $ W(-u_1,u_2) = W(u_1,u_2) $. The general existence of solutions with property similar to (A.11) for potentials with such symmetry hypothesis can be found in \cite{ABG}.

More generally, if $ a^2 +b^2 =1= c^2 +d^2 $, then
\begin{align}
u(x,y)= (\textrm{tanh}(ax+by) + \textrm{tanh}(cx+dy) , \textrm{tanh}(ax+by) - \textrm{tanh}(cx+dy))
\end{align}
solves \eqref{AllenCahnSystem} and we obtain infinitely many solutions which connect the four minima of $ W $ in sectors of variable angle.

$ \\ \\ $

\textbf{Example (3)} If $ W(u_1,u_2) = u_1^2 + u_2^2 -1 $, then
\begin{align*}
u(x,y)= ( c_1 e^{a_1x+b_1y}+ c_2 e^{a_2x+b_2y} +c_3 e^{a_3x+b_3y}+ c_4 e^{a_4x+b_4y} , 
\\ c_1 e^{a_1x+b_1y}+ c_2 e^{a_2x+b_2y} - c_3 e^{a_3x+b_3y}- c_4 e^{a_4x+b_4y})
\end{align*}
is a solution of $ \Delta u = W_u(u) $, where $ a_i^2+b_i^2 =2 \;\:,\; c_i \in \mathbb{R} $. 

In this case, $ \lbrace W = 0 \rbrace = \lbrace u_1^2+u_2^2=1 \rbrace . \\ $

Also, if $ W(u_1,u_2) = W(u_1^2+u_2^2) $ and $ W' <0 $, we have that
\begin{align*}
u(x,y) = (\textrm{cos}(ax+by+c), \textrm{sin}(ax+by+c))
\end{align*}
with $ a^2+b^2= -2W'(1) $, is a solution to $ \Delta u =W_u(u) $.
$ \\ \\ \\ $

\section{Some examples of entire solutions of the Navier-Stokes equations}

$ \\ $

First we note that some solutions of the 3D Euler equations in section 3 have the form $ u=(u_1, c_1u_1 + \tilde{c}_1,c_2u_1 + \tilde{c}_2) $, that is, we have linear dependence of the components of the solution. So, now we will determine some specific examples of solutions of the Navier-Stokes equations with linear dependent components.

Let $ u=(u_1,u_2) \;,\; u_i=u_i(x,y,t) : \mathbb{R}^2 \times (0,+ \infty) \rightarrow \mathbb{R} $ defined as
\begin{equation}\label{NSSol}
\begin{gathered}
u_1(x,y,t)= c_1 g(x-c_1y,t) -c_1A(t) +c_2 \;\;,\;\; u_2(x,y,t) = g(x-c_1y,t) -A(t) \\
\textrm{and} \;\;\; p(x,y,t) = a(t)(c_1x+y)+b(t) \;\;,t>0 \;\;\;,\; c_1,c_2 \in \mathbb{R} \\
\textrm{where} \;\;\; g_t + c_2 g_s = \mu (c_1^2+1) g_{ss} \;\;,\; g=g(s,t): \mathbb{R}^2 \rightarrow \mathbb{R} \;\; \textrm{and} \; A'(t)=a(t) \;\;,\; a,b,A: \mathbb{R} \rightarrow \mathbb{R}  
\end{gathered}
\end{equation}
then $ u $ is a solution of
\begin{equation}\label{NS2D}
\begin{cases}
u_{1 \: t} +u_1 u_{1 \: x} + u_2 u_{1 \: y} =-p_x + \mu \Delta u_1 \\
u_{2 \: t} +u_1 u_{2 \: x} + u_2 u_{2 \: y} =-p_y + \mu \Delta u_2 \\
u_{1 \: x} + u_{2 \: y} =0
\end{cases}
\;\;\; , \; \mu >0
\end{equation}
$ \\ \\ $

Similarly in the three dimensional case, we give some examples of solutions of
\begin{equation}\label{NS3D}
\begin{cases}
u_{1 \: t} +u_1 u_{1 \: x} + u_2 u_{1 \: y}+ u_3 u_{1 \: z} = - p_x + \mu \Delta u_1 \\
u_{2 \: t} +u_1 u_{2 \: x} + u_2 u_{2 \: y} + u_3 u_{2 \: z} = - p_y + \mu \Delta u_2 \\
u_{3 \: t} +u_1 u_{3 \: x} + u_2 u_{3 \: y} + u_3 u_{3 \: z} = -p_z + \mu \Delta u_3  \\
u_{1 \: x} + u_{2 \: y} + u_{3 \: z} =0
\end{cases}
\;\;\;,\; \mu >0
\end{equation}

Let $ g=g(s, \eta ,t) \;,\; g : \mathbb{R}^2 \times (0, + \infty) \rightarrow \mathbb{R} $ be a solution of
\begin{equation}
g_t - ( \frac{\tilde{c}_1}{2c_1} + \frac{\tilde{c}_2}{2c_2}) g_s +( \frac{\tilde{c}_1}{2c_1} - \frac{\tilde{c}_2}{2c_2}) g_{\eta} = \mu (\frac{1}{4c_1^2} + \frac{1}{4c_2^2}) (g_{ss} +g_{\eta \eta}) + \mu g_{ss}
\end{equation}
where $ \mu >0 \;,\; c_1,c_2, \tilde{c}_1 , \tilde{c}_2 \in \mathbb{R} $ and $ t>0 $.

Then $ u=(u_1,u_2,u_3) \;\;,\; u_i : \mathbb{R}^3 \times (0,+ \infty) \rightarrow \mathbb{R} \;\;,\; i =1,2,3 $ defined as
\begin{equation}\label{NS3DSol}
\begin{gathered}
u_1(x,y,z,t) = g( x - \frac{c_2y + c_1z}{2c_1c_2}, \frac{c_2y - c_1z}{2c_1c_2},t) - A(t) \;\;\;,\; (x,y,z) \in \mathbb{R}^3 \;,\; t>0 \\
u_2 (x,y,z,t) = c_1 u_1(x,y,z,t) + \tilde{c}_1 \;\;\;,\;\; u_3(x,y,z,t) = c_2 u_1(x,y,z,t) + \tilde{c}_2 \\
\textrm{and} \;\; p(x,y,z,t) = a(t)(x + c_1y +c_2z) +b(t) \\
\textrm{where} \;\; A'(t) =a(t) \;,\; a,A : \mathbb{R} \rightarrow \mathbb{R}
\end{gathered} 
\end{equation}
is a solution of \eqref{NS3D}. $ \\ $

Therefore we conclude to the following $ \\ $

\begin{proposition}\label{NSlinearsolutions} Let $ u=(u_1,u_2,u_3) \;\:,\; u_i,\: p: \mathbb{R}^3 \times (0, + \infty) \rightarrow \mathbb{R}^3 $ and consider the initial value problem
\begin{equation}\label{NS3DIV}
\begin{cases}
u_t +u \nabla u = - \nabla p + \mu \Delta u \\
\textrm{div} \: u =0 \\
\lim_{t \rightarrow 0^+} u(x,y,z,t) = h(x,y,z)
\end{cases}
\;\;\;,\; \mu >0 \;\;,\; (x,y,z,t) \in \mathbb{R}^3 \times (0, + \infty)
\end{equation}
where $ h=(h_1,c_1h_1 + \tilde{c}_1, c_2h_1 +\tilde{c}_2) $ and $ h_1(x,y,z)= H(2c_1c_2 x - c_2y - c_1z) \;\;,\; c_1,c_2, \tilde{c}_1 , \tilde{c}_2 \in \mathbb{R} $ such that $ \tilde{c}_1 c_2 + c_1 \tilde{c}_2 =0 $ and $ H $ smooth.

Then there exists a smooth, globally defined in $ t>0 $, solution to \eqref{NS3DIV}. 

In particular,
\begin{equation}\label{NS3DIVSol}
\begin{gathered}
u(x,y,z,t) = (u_1,c_1u_1 + \tilde{c}_1, c_2u_1+ \tilde{c}_2) \;\: \textrm{and} \;\: p(x,y,z,t)=a(t)(x+c_1y+c_2z) +b(t) \\
\textrm{where} \;\; u_1(x,y,z,t)= g(2c_1c_2x - c_2y -c_1z,t) -A(t) \\ \textrm{and} \;\; g=g(s,t)= \frac{1}{2\sqrt{\pi t}} \int_{\mathbb{R}} e^{-\frac{|s-w|^2}{4 \tilde{\mu}t}} H(w)dw \;\;,\; \tilde{\mu}= \mu (4c_1^2 c_2^2 +c_1^2 +c_2^2) \\ ( A'(t)=a(t) \;\: ,\; A(0)=0)
\end{gathered}
\end{equation}
\end{proposition}

$ \\ $

\textbf{Remark B.1} We can also have the same result for a bit more general initial values $ h $ in Proposition 3, as we can see from (B.4),(B.5). It suffices to have linear dependency of the components of $ h $ and $ h_1 $ above can also be for example of the form $ h_1(x,y,z) = H(2c_1c_2x-c_2y-c_1z,c_2y-c_1z) $.

\end{appendix}

$ \\ $

\end{document}